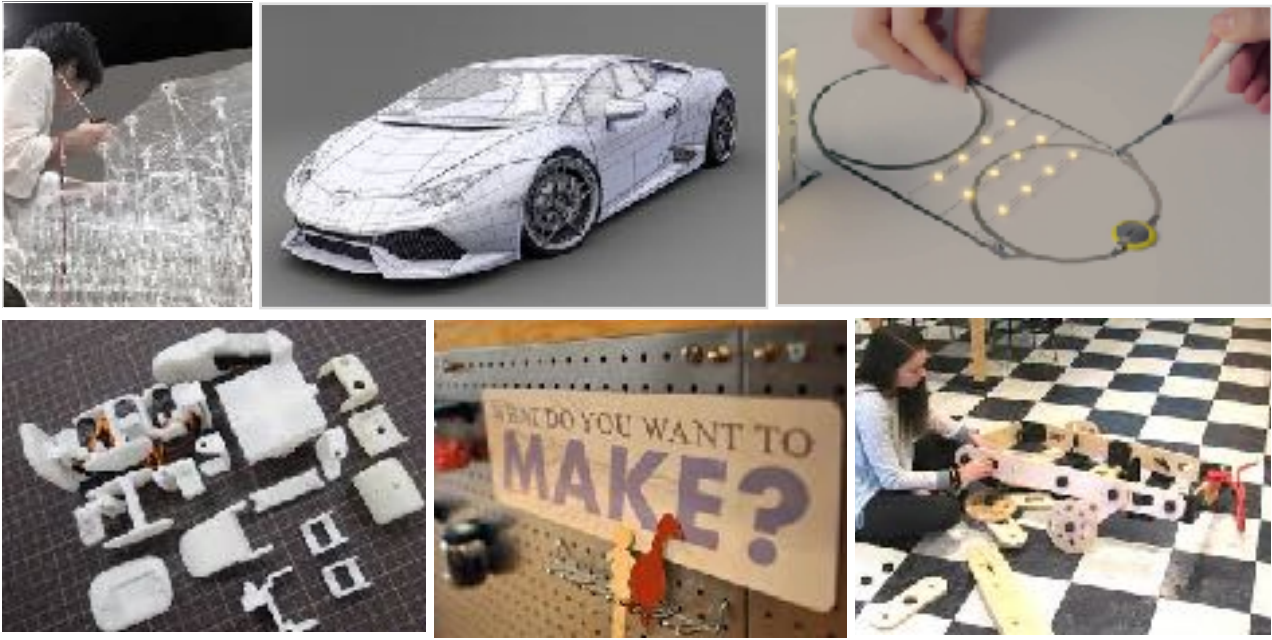

Implementación de un

# STEAM-Makerspace

para dinamizar los procesos de enseñanza
y aprendizaje de las

# matemáticas
(con énfasis en la geometría)

de los estudiantes de nivel medio superior

———

Línea de investigación:

**Aprendizaje y logro educativo de las ciencias y tecnología**

Luis Adrían Martínez Pérez[1]

———

[1] Responsable Técnico, email: lmartinez@edu.prp.ceyca.com, lamp@comunidad.unam.mx

**PROYECTO**: Implementación de un STEAM *Makerspace* para dinamizar los procesos de enseñanza y aprendizaje de la matemática (con énfasis en la geometría) de los estudiantes de segundo año del nivel medio superior.

**1. LÍNEA DE INVESTIGACIÓN**

Aprendizaje y logro educativo de las ciencias y tecnología.

**2. ANTECEDENTES TEÓRICOS**

> «Ἀγεωμέτρητος μηδείς εἰσίτω» (*'No entre nadie que no sepa geometría'*). **Platón, 388 a.C.**
>
> *Nos encontramos con patrones todo el tiempo, diariamente: en el lenguaje hablado y escrito, en las formas musicales y las imágenes de vídeo, en el diseño ornamental y la geometría natural, en las normas de vialidad y en los objetos que construimos. Nuestra habilidad para reconocer, interpretar y crear patrones constituye un elemento clave para tratar con el mundo que nos rodea.* **Senechal, M.**[2]

Tal como señalan las anteriores citas, la estricta importancia que tiene el desarrollo de las competencias matemáticas no sólo para el avance del pensamiento científico, el saber humanista o la expresión artística y creativa, sino para la vida diaria, es innegable.[3] Tomando en cuenta lo anterior y a la luz de los resultados de las pruebas nacionales e internaciones más recientes[4], sabemos que nos encontramos frente a una problemática educativa nacional: las deficiencias que presentan los jóvenes que concluyen la educación media superior, en los conocimientos y habilidades matemáticas, nos revela que se afectará, de forma negativa, el desarrollo intelectual, profesional y personal de los mismos a futuro.[5] El Proyecto Educativo del Área de Ciencias del nivel medio superior, en el que un STEAM-*Makerspace* tiene importante cabida, va encaminado a entender cómo subsanar dicha deficiencia.[6]

---

[2] Senechal, M., *Forma. La enseñanza agradable de las matemáticas* pp. 149. Ed. Limusa. México. 2004.
[3] Dicha relación fundamental ya ha sido ampliamente documentada. Para la relación arte-matemáticas véase Edo B., Mercè (2005, 2009); Dos Santos, J. Manuel (2009); Giménez, Joaquín (2009); Kaplan, A., et.al., (2015). Para la relación matemáticas-ciencias véase Carrillo M., y Sánchez, Antonio (2006); Rodriguez, Elena M. (2011); Trejo, Elia et.al. (2013); (2016). Para la relación matemáticas-humanidades ver . Sobre las matemáticas en la vida diaria ver Alsina A., Callis, J., et.al. (1998); Vicente, S., Dooren, W. y Verschaffel, L. (2008); PISA (2012); Solares, D., Solares, A. y Padilla, Erica (2016).
[4] Para consultar los resultados del examen PLANEA (2015-16) (2016-17) véase http://planea.sep.gob.mx/ms/estadisticas/resultadoscontextualizados/infoNacional.php?campo=null; para consultar los resultados de la prueba PISA (2015) (2016) verse https://www.oecd.org/pisa/PISA-2015-Mexico-ESP.pdf . En el caso del Colegio Ceyca, la Coordinación de Ciencias ha realizado los análisis de PLANEA y College Board (2015) (2016) en éste rubro. Para una retrospectiva de la educación matemática en México, véase Ávila, A (2016)
[5] Para la relación matemáticas-desarrollo cognitivo, creativo y productivo, véase a Rosskopf., Myron J., (1971); Wagner, S. (1981); Collis, K., (1981); Strauss, S. (2000); Gurganus, S.P. (2007); Arch Tirado (2008); Vargas, G. (2011); PISA (2012); (2014); Siegler, R. S. (2016).
[6] Dicho proyecto, que en la actualidad está en proceso de implementación, es resultado del análisis, durante los dos últimos ciclos (2015-2016 y 2016-2017), de los resultados de las pruebas nacionales, internacionales y de la propia red de colegios, así como de la reflexión de la experiencia docente en ésta área. Nuestro Proyecto Educativo enfatiza una respuesta radical al problema señalado, una respuesta que, desde luego, implica la incorporación de las tecnologías educativas de vanguardia y que aquí, en el Colegio Ceyca y en la Coordinación de Ciencias, hemos insistido en impulsar.

Las causas documentadas del porqué de ésta problemática son diversas.[7] Sin embargo, una nos resulta fundamental: con las reformas realizadas al curriculum de matemáticas en las décadas de 1960 a 1980 se registró una disminución del rol de la geometría en la formación educativa.[8] A nuestro juicio, apoyado en lo documentado[9] y como resultado de nuestra experiencia docente[10], consideramos que este soslayo ha derivado en una deficiente comprensión de las matemáticas y, en consecuencia, en un bajo rendimiento académico de los estudiantes no sólo en esa área, sino en todas las ciencias que dependen de ellas directa o indirectamente.[11]

Recuperar la geometría en el currículo escolar y volverla una experiencia clara, evidente, interesante y divertida para los jóvenes, a través de proyectos interdisciplinarios realizados por medio de los recursos tecnológicos de vanguardia -como los que ofrece un STEAM *Makerspace*-, podría resultar fundamental para resarcir el rezago matemático en nuestros jóvenes[12] (ESTA ES LA TESIS FUNDAMENTAL DEL PROYECTO EDUCATIVO GENERAL, EN DONDE EL MAKER SPACE APARECE COMO UNA HERRAMIENTA DE APOYO). A nuestro juicio, implementar un STEAM *Makerspace* para dinamizar los procesos de enseñanza y aprendizaje de las matemáticas es reconocer, con seriedad, que el ambiente tecnológico en el que se desarrolla la vida y el futuro de nuestros jóvenes exige que los centros educativos y los docentes incorporen estas mismas herramientas a su infraestructura y prácticas respectivamente.[13]

### 3. Objetivo general del proyecto y resultados esperados

El objetivo general es que los estudiantes del segundo año de bachillerato, empleen las herramientas de un STEAM-*Makerspace* para dinamizar sus procesos de enseñanza y aprendizaje de las matemáticas, específicamente de la geometría.

Los resultados que esperamos son:

a) que el alumno observe, modele, visualice y simule matemáticamente, a partir de las

---

[7] Véase: Planas, N. (2016);

[8] Véase: ICEM, (2001);

[9] La geometría está siendo fuertemente defendida, en los más diversos países como método de enseñanza de las matemáticas en todos los niveles de escolaridad (Salett B., M. y Nelson, H (2004); Corica, A.R. y Marin, E.A. (2014);        ). Véase también para la importancia de la geometría en la formación educativa a Santaló, Luis A. (1994); Barrera, F. Y Santos, M. (2001); Edo (2005); Goncalves, R. (2006); Bogomolny (2010); Vargas, G. y Gamboa, R. (2012); Ospitaletche-Borgmann y Martínez (2012); Arnal, A. y Planas, N. (2013); Acosta, M. y Camargo, L. (2014); Rojano, T. (2014). Sobre fundamentos teóricos de la geometría y su importancia en el edificio matemático, véase Hilbert (1993); Husserl (2000); Corre (2006); Giovannini (2014).

[10] Como se ha mencionado arriba, el Proyecto Educativo de la Coordinación de Ciencias en el que se enmarca este proyecto, engloba las reflexiones sobre nuestra experiencia docente en la enseñanza de las matemáticas, específicamente de la geometría, así como los análisis que hemos realizados.

[11] Como afirmó M. Kline: "*En el lugar de las matemáticas hay muchas moradas, y entre ellas, la más elegante es la Geometría*". Pies de página arriba, hemos hecho referencia a los estudios sobre la relación e impacto de las matemáticas con las diversas esferas del pensamiento y quehacer humano. Para la relación matemáticas, ciencias y rendimiento escolar ver        . También puede revisarse el Análisis (2016) realizado por el área de Ciencias en relación al rendimiento de las materias científicas el Ceyca.

[12] Para los estudios sobre el impacto del uso de la tecnología en la enseñanza matemática, véase: Williamson, S. y Kaput, J. (1999); Martín, W (2000); Lupiánez, J. y Moreno, L. (2001); Sandoval, I. T. y Moreno, L.S. (2012); Parra, O. y Díaz, V. (2014); Rodriguez, S., Ibáñez, J., Valencia, N. (2016).

[13] Como hemos señalado anteriormente, el Proyecto Educativo de la Coordinación de Ciencias del Colegio Ceyca tiene presente la importancia del ambiente tecnológico presente y futuro de los estudiantes.

herramientas de un STEAM-*Makerspace*, fenómenos naturales y sociales.

b) que el alumno incorpore los modelos, visualizaciones, simulaciones matemáticas desarrolladas en las herramientas de un STEAM-*makerspace* en sus proyectos de investigación científica.

c) que el alumno exprese de forma creativa, por medio de las herramientas de un STEAM-*Makerspace*, la relación entre matemáticas, arte y ciencias.

d) que el aprendizaje o corrección de los conceptos matemáticos a través de estos proyectos, posibiliten al alumno elevar su interés, gusto y rendimiento académico por las ciencias.

e) que el docente incorpore las tecnologías, metodologías y dinámicas *maker* para los procesos de enseñanza y aprendizaje.

**4. METODOLOGÍA**

Para implementar un STEAM-*makerspace* en el proceso de enseñanza y aprendizaje, elegiremos, bajo parámetros definidos, un <grupo piloto> de segundo año de bachillerato, durante dos años. Los profesores del área de ciencias y tecnología diseñarán una serie de proyectos interdisciplinarios (con énfasis en el aprendizaje de los conceptos geométricos) que los alumnos desarrollarán de manera activa, creativa y lúdica. Estos proyectos se diseñarán bajo un conjunto parámetros e indicadores para que, posteriormente y mediante el Ciclo de Modelación de Rodríguez[14], podamos analizar y evaluar el resultado TANTO DE LA MEJORA DE LA COMPRENSIÓN DE LA GEOMETRÍA COMO de la implementación de un STEAM-*makerspace* en la dinamización de los procesos de enseñanza y aprendizaje de las matemáticas (con énfasis) en geometría, tanto para alumnos como para profesores.

**5. PRESUPUESTO**

El precio por *servicios* de asesoría de un STEAM-*makerspace*, de la empresa Hacedores Makerspace, que incluye el diseño, planeación y recomendaciones del *makerspace* es de $75,000.00 más IVA. En este rubro recomiendan agregar un presupuesto adicional para adaptación de imagen o *look and feel* del espacio como proyecto arquitectónico de $150,000.00 más IVA.

La inversión estimada para el Equipamiento del *makerspace*, esto es herramientas, equipos, materiales y mobiliario, dependerá de la selección que se realice en la primera fase, esta cantidad va en una estimación de los $450,000 pesos hasta los

---

[14] Proponemos el Ciclo de Moderación de Rodriquez porque en éste "*…se analiza el proceso de moderación matemática como un conjunto de diversas praxeologías. Esta idea se retoma de la teoría antropológica de lo didáctico propuesta por Yves Chevallard en 1985, cuyo principal postulado establece la idea de que toda actividad humana puede describirse a través de praxeologías. Una praxeología de acuerdo con Chevallard alberga cuatro elementos principales: tarea, técnica, tecnología y teórico. Los dos primeros conforman el bloque práctico-técnico, relativo al saber hacer, y los dos segundos, referentes al saber, conforman el bloque tecnológico-teórico*" pp. 104

$600,000 dependiendo de la cantidad y calidad de equipamiento que el Colegio defina con nuestra ayuda.

El costo de la Formación de Capital Humano (tanto del Operador del *makerspace*, como del grupo de trabajo o Docentes):

a) Consultoría en la selección y capacitación del operador será de $20,000.00 más IVA. Incluye: 20 horas de Capacitación de operación del *makerspace* para el operario.

b) Programa de Formación Docente: sesiones de Inducción al Movimiento *Maker* y sesiones prácticas de cómo diseñar Dinámicas *Maker*: $10,000.00 + IVA x sesión de 6 horas por grupo de 5 profesores.

Acompañamiento personalizado en atención de vinculación con otras instituciones, comunicación institucional y promoción externa en manejo de mensajes sobre la nueva cultura *maker* durante el primer periodo de 3 meses: $25,000.00 más IVA.

Cursos Académicos Maker para alumnos diseñados por Hacedores, a impartir durante un ciclo escolar. Ejemplo, *Music Makers*, *MeteorITo* o De Cero a *Maker*: $14,000.00 más IVA al mes por grupo de máximo 25 alumnos impartido en 35 clases a lo largo del año escolar.

**Total de inversión: $944,000.00** (*Este total NO incluye IVA y no contemplan gastos de viaje ni costos de fletes en el caso de equipamiento y serán calculados adicionales a los servicios, por lo que el total podría ser superior a la cantidad señalada. En todo caso, depende de los acuerdos con la empresa Hacedores *Makerspace* )

### 6. Impacto y mecanismo de transferencia

Para diseñar las mejores estrategias de transferencia, primero proponernos integrar una memoria de trabajo. Posteriormente proponemos que los resultados del análisis de las memorias se difundan a toda la red a través de medios digitales, mesas de trabajo, seminarios y conferencias. A la luz de los resultados publicados y de los debates generados la red, o cada colegio, evaluará la pertinencia de implementar un STEAM *makerspace*. En ese caso, se propone que se defina, en conjunto con los colegios interesados, el mecanismo idóneo para la transferencia de los conocimientos y experiencias obtenidas.

### 7. Bibliografía citada